\documentclass[12pt]{article}
\usepackage[latin2]{inputenc}
\usepackage{amsfonts, amsmath, amsthm, mathtools, color}
\usepackage{enumerate}
\usepackage{amssymb}
\usepackage{authblk}
\usepackage{thmtools} % new by Moritz
\usepackage[T1]{fontenc}
\usepackage{bbm}%for\mathbbm{1}

\usepackage[colorlinks=true, urlcolor=blue, linkcolor=blue,
citecolor=black]{hyperref}
\usepackage{color}
\usepackage[english]{babel}
  
\addto\extrasenglish{%

}

\usepackage[backgroundcolor=white, bordercolor=blue, linecolor=blue]{todonotes}
\renewcommand{\d}{\mathrm{d}}

\newcommand{\dx} {\, \mathrm{d}x}
\newcommand{\dy} {\, \mathrm{d}y}
\newcommand{\dz} {\, \mathrm{d}z}
\newcommand{\dw} {\, \mathrm{d}w}
\newcommand{\dist} {\mathrm{dist}}
\newcommand{\LL} {\mathrm{L}}

\theoremstyle{plain}
\declaretheorem[title=Theorem]{theorem}

\declaretheorem[title=Corollary,sibling=theorem]{corollary}

\def\XXint#1#2#3{{\setbox0=\hbox{$#1{#2#3}{\int}$ }
\vcenter{\hbox{$#2#3$ }}\kern-.6\wd0}}

\theoremstyle{definition}
\declaretheorem[title=Definition,sibling=theorem]{definition}
\declaretheorem[title=Remark,sibling=theorem]{remark}
\declaretheorem[title=Remark, numbered=no]{remark*}
\declaretheorem[title=Example, sibling=theorem]{example}

\numberwithin{equation}{section}

\usepackage{dsfont}
\newcommand{\vis}{{\operatorname{vis}}}
\newcommand{\cen}{{\operatorname{cen}}}
\newcommand{\cE}{{\ensuremath{\mathcal{E}}}}
\newcommand{\cF}{{\ensuremath{\mathcal{F}}}}

\newcommand{\cW}{{\ensuremath{\mathcal{W}}}}

\newcommand{\R}{{\ensuremath{\mathbb{R}}}}
\newcommand{\N}{{\ensuremath{\mathbb{N}}}}

\renewcommand{\P}{\ensuremath{\mathbb{P}}}
\renewcommand{\dj}{d\kern-0.4em\char"16\kern-0.1em}

\renewcommand{\proof}{\noindent\textbf{Proof.}\ } %changed by Moritz

\newcommand{\proofof}{\noindent\textbf{Proof of}\ }

\renewcommand{\qed}{\hfill\ensuremath{\Box}\\} % changed by Moritz

\newcommand{\diam}[1]{\mathop{\textrm{diam}} #1}
\renewcommand{\dist}[1]{\mathop{\textrm{dist}} #1}

\usepackage{graphicx}
\usepackage{subcaption}
\usepackage{wrapfig}
\usepackage{xfrac}
\allowdisplaybreaks

\title{Nonlocal quadratic forms with visibility constraint}

\author{{Moritz Kassmann}\thanks{The authors gratefully acknowledge financial support by the German Science Foundation DFG via SFB1283.}\, \thanks{Universit\"{a}t Bielefeld, Fakult\"{a}t f\"{u}r Mathematik, Postfach 100131, D-33501 Bielefeld, Germany, \emph{moritz.kassmann@uni-bielefeld.de}}
\quad and
\quad {Vanja Wagner$^\ast$\thanks{University of Zagreb, Faculty of Science, Department of Mathematics, Bijeni\v cka cesta 30, 10000 Zagreb, Croatia, \emph{wagner@math.hr}}}
}

\begin{document}

\maketitle

\abstract{Given a subset $D$ of the Euclidean space, we study nonlocal 
quadratic forms that take into account tuples $(x,y) \in D \times D$ if and 
only 
if the line segment between $x$ and $y$ is contained in $D$. We discuss 
regularity of the corresponding Dirichlet form leading to the existence of a 
jump process with visibility constraint. Our main aim is to investigate 
corresponding 
Poincar\'{e} inequalities and their scaling properties. For dumbbell shaped domains we show that the forms satisfy a Poincar\'{e} 
inequality with diffusive scaling. This relates to the rate of convergence of eigenvalues in singularly perturbed domains.}

\medskip

\noindent{\it AMS Subject Classification (MSC2010)}: 46E35, 47G20, 60J75, 26D15 \\
\noindent{\it Keywords}: Function spaces, nonlocal Dirichlet forms, Poincar\'{e} inequalities, Integro-differential operators, jump processes

\pagebreak[3]

\section{Introduction}\label{sec:intro}
\subsection{Motivation and Setup}
The aim of this work is to study nonlocal quadratic forms related to Markov jump processes, corresponding function spaces, and Poincar\'{e} inequalities. Let us begin with a simple example. If $D = D^- \cup D^+ \subset \R^d$ is the union of two disjoint components, then any diffusion on $D$ decomposes into two separate diffusions. For a jump process, say an isotropic $\alpha$-stable process, this is different because the connectedness of the domain is irrelevant for the jump process. If  $D = D^- \cup \Gamma \cup D^+$, where $\Gamma$ is a thin corridor connecting the two components $D^-, D^+$, then a diffusion has to pass through $\Gamma$ in order to reach one component from the other. This has led to interesting quantitative studies of eigenvalue problems for generators of diffusions in dumbbell shaped domains. Very similar situations appear in the study of meta-stability when a diffusion has to overcome a hill in order to move from one well of the considered energy landscape to another one. 

\medskip

Similar problems for generators of jump processes seem to be uninteresting because the jump process does not need to pass through the thin corridor in order to move from $D^-$ to $D^+$. In this work we introduce and study nonlocal quadratic forms that generate jump processes that do have this property. Jumps between points $x \in D$ and $y \in D$ can only take place if the line segment between $x$ and $y$ is contained in the domain $D$, i.e., if the points are "visible" one from another. Our focus is on Poincar\'{e} inequalities in non-convex domains of the form $D = D^- \cup \Gamma \cup D^+$, i.e., so called dumbbell shaped domains.

\begin{figure}[h]\hypertarget{figure1}{ }
	\centering
\includegraphics[width=0.25\textwidth]{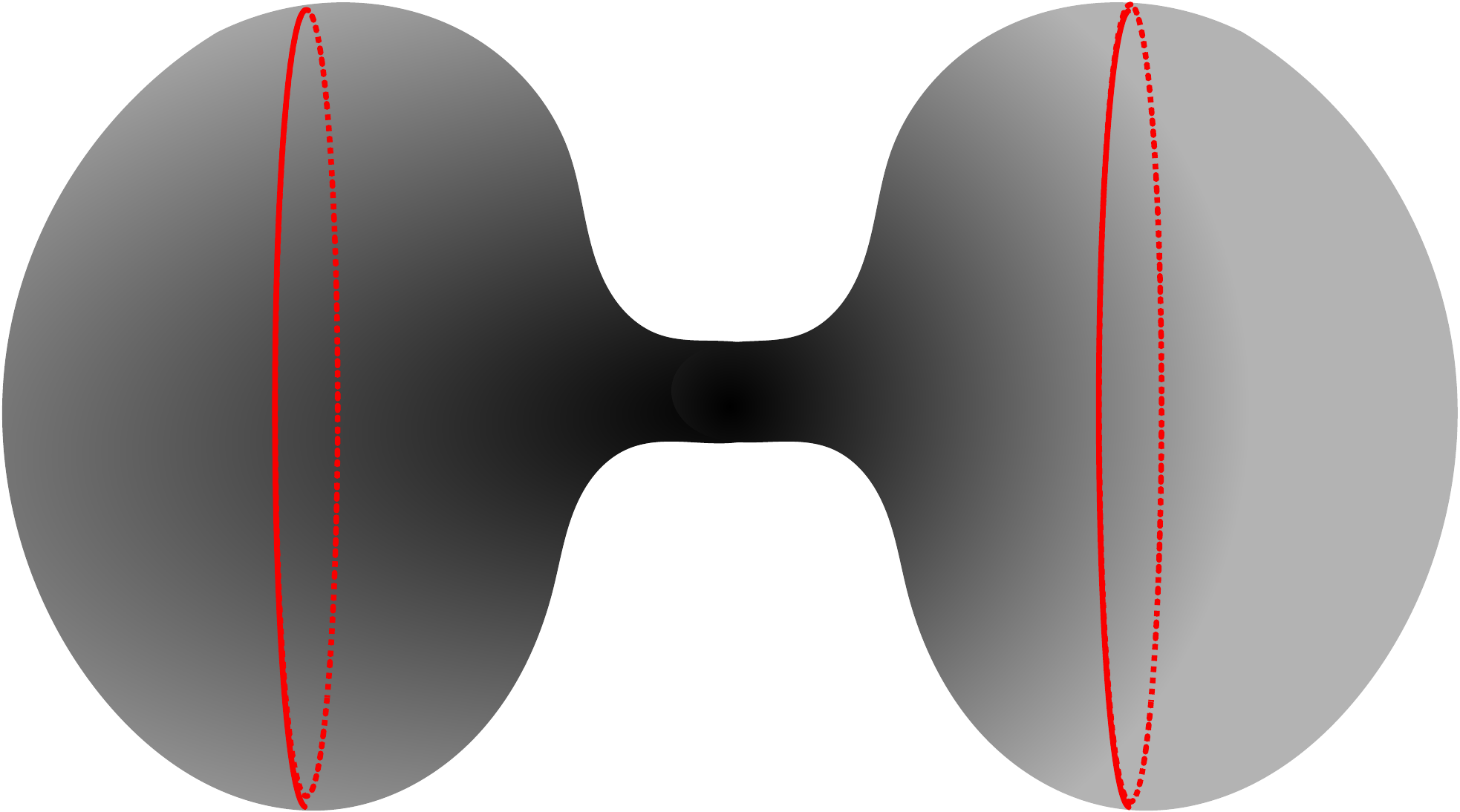}
\caption{A dumbbell shaped domain}
\end{figure} 

We consider sequences of such domains, where the corridor $\Gamma$ is fixed  and the sets $D^-$, $D^+$ are assumed to be growing. In \autoref{thm:poincare-intro} we establish Poincar\'{e} inequalities in such domains. It turns out that the scaling behavior of the Poincar\'{e} constant is identical for local diffusive-type quadratic forms and nonlocal jump-type quadratic forms, no matter the value of $\alpha \in (1,2)$. This phenomenon shows that the visibility constraint has a serious impact. The probabilistic interpretation behind this phenomenon is that the visibility-constrained jump process is forced to pass through the corridor in order to move from one component to the other one. On large scales, this restriction makes the jump process as slow as the Brownian Motion, compare \autoref{thm:poincare-local} and \autoref{thm:poincare-nonlocal}. For a special class of domains [cf. \eqref{eq:Condition_C}] and $\alpha \in (0,1)$ we show that the long-range connections of nonlocal forms may lead to a different scaling, see the second part of \autoref{thm:poincare-intro}. 

\medskip

\pagebreak[3]

The study of Poincar\'e inequalities for large domains of the aforementioned type is closely connected with the study of eigenvalues in bounded singularly perturbed domains. Here one considers a sequence of domains $\Omega_\epsilon$ together with some limit domain $\Omega_0 \subset \Omega_\epsilon$ where the Lebesgue measure of $\Omega_\epsilon \setminus \Omega_0$ tends to zero as $\epsilon$ tends to zero. The study of eigenvalue problems with Neumann or Dirichlet data in such domains has a long history, see  \cite{Ars71}, \cite{Ars72}, \cite{Bea73} \cite{LoSa79},  \cite{Jim89} and \cite{Arr95}. Related problems concern Helmholtz resonators, see \cite{HiMa91}, \cite{BHM94}, \cite{BHM95}. Such problems do not make any sense when one considers eigenvalue problems with respect to classical nonlocal operators like the fractional Laplace operator. This is because the nonlocal operator does not at all react to, resp.~"feel", the singular perturbation. One motivation for the introduction of nonlocal operators with visibility constraint is to  study such problems also in the framework of nonlocal operators.

\medskip

Let us set up the mathematical context. Throughout this paper we assume that $D \subset \R^d$ is a measurable subset and $k: D \times D \setminus \operatorname{diag} \to [0, 
\infty)$ is a measurable function such that
\begin{align}\label{eq:integrability1}
\sup_{x\in D}\int_{D\setminus\{x\}}(1\wedge|x-y|^2)k(x,y)\dy<\infty.
\end{align}
We consider Hilbert spaces of the form 
$H(D) = \big\{f \in L^2({D}) \big|\, |f|_{H(D)} < \infty \big\}$ 
with a semi-norm $|f|_{H(D)}$ given by 
\begin{align}\label{eq:seminorm-H}
 |f|^2_{H(D)} = \iint\limits_{D D} (f(y)-f(x))^2 
k(x,y) \dx \dy  \,.
\end{align}
Note that the condition \eqref{eq:integrability1} on $k$ implies that $C_c^\infty(D)\subset H(D)$. The space $H(D)$ is endowed with the norm $\|f\|_{H(D)}$ given by 
$\|f\|_{H(D)}^2 = \|f\|_{L^2(D)}^2 + |f|^2_{H(D)}$. Note that, without loss 
of generality, one can assume the function $k$ to be symmetric due to the 
symmetric structure of the double-integral in \eqref{eq:seminorm-H}.

\medskip

Let us look at some 
special choices of $k(x,y)$. If $k$ is a bounded 
function then $H(D)$ equals $L^2(D)$. If $k(x,y) = c_{s, d} 
|x-y|^{-d-2s}$, for some $s \in (0,1)$ and some positive $c_{s, 
d}$, then $H(D)$ equals the well-known Sobolev-Slobodeckii space 
$H^{s}(D)$. Note that, when relating $k(x,y)$ resp.~the function space $H^{s}(D)$ to stochastic processes, it is common to replace $2s$ by $\alpha \in (0,2)$ because the corresponding jump process is called $\alpha$-stable jump process. Since this work is mostly concerned with quadratic forms and functional inequalities, we use $s \in (0,1)$. One can choose the constant $c_{s, 
d}$ such that for every smooth function $f$ the norm $\|f\|_{H^{s}(D)}$ 
converges to $\|f\|_{L^2(D)}$ as 
$s \to 0+$ and $\|f\|_{H^{s}(D)}$ converges to $\|f\|_{H^1(D)}$ as 
$s \to 1-$. Here $H^1(D)$ equals the classical Sobolev space of 
all functions $f \in L^2({D})$ with generalized derivative $\nabla f \in 
L^2(D;\R^d)$.

\medskip

This work is concerned with a new kind of nonlocal spaces. For $x\in D$ let 
$D_x\subset D$ be the visible region in $D$ from point $x$, i.e.
\begin{align*}
 D_x=\{y\in D|\, tx+(1-t)y\in D,\,\forall t\in(0,1)\}.
\end{align*}
Given a symmetric function $k: D \times D \setminus \operatorname{diag} \to [0, 
\infty)$ as before, we introduce a smaller semi-norm on $L^2(D)$ induced by the bilinear form 
\begin{align}\label{eq:seminorm-Hvisible}
\cE^{\operatorname{vis}}(f,f) &:= \iint\limits_{D D_x} (f(y)-f(x))^2 
k(x,y) \dx \dy \\
&= \frac12 \iint\limits_{D D} (f(y)-f(x))^2 
k(x,y) \big( \mathds{1}_{D_x}(y) +  \mathds{1}_{D_y}(x) \big) \dx \dy \,,\nonumber
\end{align}
if this quantity is finite. 
\begin{figure}
	\centering
	\begin{subfigure}[b]{0.45\textwidth}
		\hspace*{4ex}\includegraphics[width=0.75\textwidth]{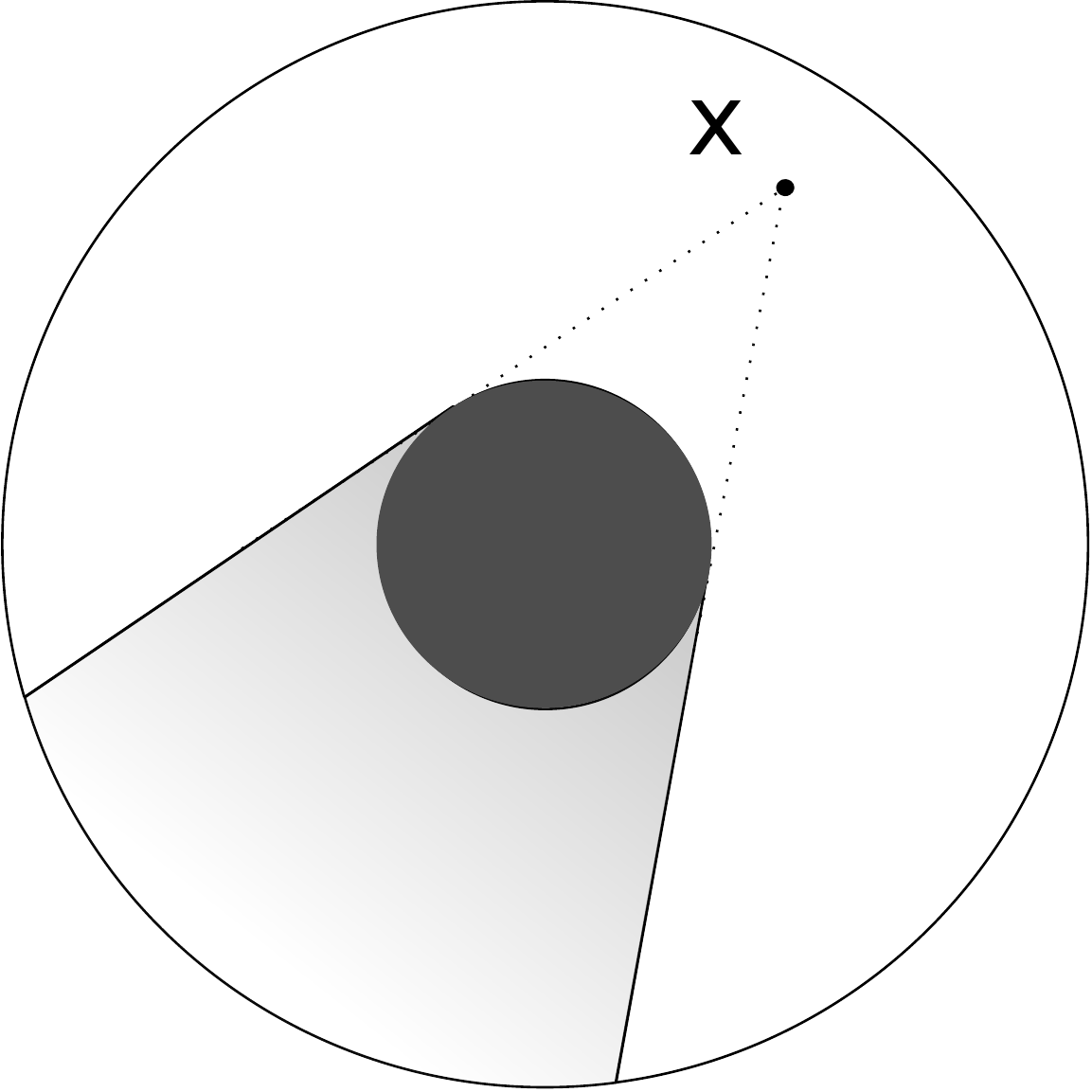}
		\caption{An annulus: $B_1 \setminus B_{1/3}$}
		\label{fig:gull}
	\end{subfigure}
	~ %add desired spacing between images, e. g. ~, \quad, \qquad, \hfill etc. 
	\begin{subfigure}[b]{0.45\textwidth}
		\hspace*{5ex}\includegraphics[width=0.7\textwidth]{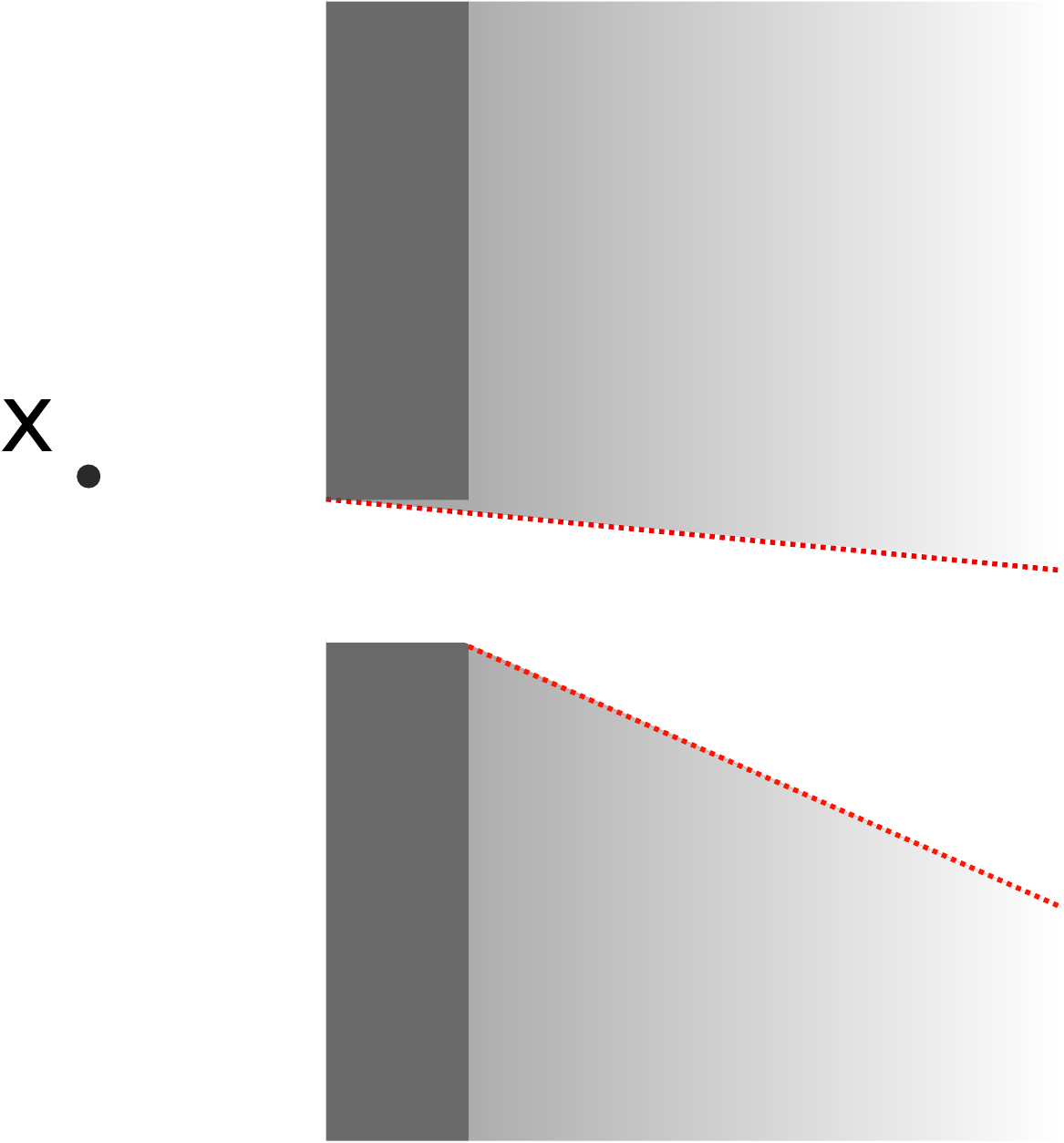}
		\caption{An unbounded domain as in \autoref{ex:1}}
		\label{fig:tiger}
	\end{subfigure}
	\caption{Two non-convex sets $D$ with a point $x$ and its region of visibility}\label{fig:domains}
\end{figure}
In these semi-norms, tuples $(x,y) \in D \times D$ are considered only 
if $x$ is an element of $D_y$ or, equivalently, $y$ is an element of $D_x$. One can imagine 
points to be connected only if they can ''see`` each other. In this sense, we 
decide to call the object defined in \eqref{eq:seminorm-Hvisible} a 
\emph{quadratic nonlocal form with visibility constraint}. Obviously, for 
convex domains $D$, this semi-norm is equal to the one defined in 
\eqref{eq:seminorm-H}, i.e.~in this case
\begin{align}\label{eq:Ecensored}
\cE^{\vis}(f,f) = \cE^{\cen}(f,f) &:= 
\iint\limits_{D D} (f(y)-f(x))^2 
k(x,y) \dx \dy \,. 
\end{align}
Here ''cen`` stands for ''censored`` and indicates that points outside of $D$ 
are not taken into account. Censored forms and corresponding stochastic 
processes have been introduced in \cite{BBC03}. Given any bilinear form $\cE$ on $L^2(D) \times L^2(D)$ as above, we set $\cE_1(f,f) = \cE(f,f) + 
\|f\|_{L^2(D)}$ as usually done. We can now define the function spaces that are of 
particular interest to us.

\begin{definition}
Let $D \subset \R^d$ be an open measurable subset of $\R^d$ and $k: D \times D 
\setminus \operatorname{diag} \to [0, \infty)$ be a measurable function. Then we 
define four function spaces:
\begin{alignat*}{2}
\cF^{\vis} &=\overline{C_c^\infty(D)}^{\cE^{\vis}_1}, \qquad &
\widetilde{\cF}^{\vis} &= \{ f \in L^2(D) |\,  \cE^{\vis}(f) < \infty \} \,, \\
\cF^{\cen} &=\overline{C_c^\infty(D)}^{\cE^{\cen}_1}, \qquad &
\widetilde{\cF}^{\cen} &= \{ f \in L^2(D) |\,  \cE^{\cen}(f) < \infty \} 
\,.
\end{alignat*}
\end{definition}

Note that the choice of the set $D$ and the kernel $k$ is very important for 
these domains. The model case that we have in mind is given by a bounded 
open non-convex set $D$ with a smooth boundary and $k(x,y) = |x-y|^{-d-2s}$ for 
some $s \in (0,1)$. We will allow for more general domains and for more general 
kernels, including weakly singular, but the main new results like \autoref{thm:poincare-intro} are new even in this model case.  

\medskip
 
One driving idea behind this project is the connection of Dirichlet forms to Markov jump 
processes. Let us recall that the pair $(\cE,\cF)$ is called a Dirichlet form on $L^2(E)$, for $E\subset \R^d$ open or closed, if $\cF$ is a dense linear subspace of $L^2(E)$ and $\cE$ is a bilinear symmetric closed form on $\cF\times\cF$ which is also Markovian, e.g., if for every $u \in  \cF$ the function $v= (u \wedge 1) \vee 0$ belongs to $\cF$ and satisfies $\cE(v,v) \leq \cE(u,u)$. See \cite[Section 1.1]{FOT94} for related definitions and examples. 
A Dirichlet form $(\cE, \cF)$ on $L^2(E)$ is 
called regular if $C_c(E) \cap \cF$ is dense in 
$C_c(E)$ w.r.t.~the supremum norm as well as in $\cF$ 
w.r.t.~the norm $\cE_1(u,u)^{1/2}$. A major result due to M.~Fukushima is 
that every regular Dirichlet form $(\cE, 
\cF)$ on $L^2(E)$ corresponds to a symmetric 
strong Markov process on $(E, \mathcal{B}(E))$, cf. \cite[Theorem 
7.2.1]{FOT94}. Note that the 
rotationally symmetric $2s$-stable L\'{e}vy process, $s\in(0,1)$, is the strong Markov process 
that corresponds to the regular Dirichlet form on $L^2(\R^d)$ defined by 
\begin{align*}
(f,g)
  \,\,\mapsto\,\,
  \iint\limits_{\R^d \R^d}
  \left(f(y)-f(x) \right) \left(g(y)-g(x) \right) \; |x-y|^{-d-2s} \dx \dy 
  \qquad (f,g\in H^s(\R^d))\,.
\end{align*} 
The tuple $(\cE^{\vis},\cF^{\vis})$ is by construction a regular 
Dirichlet form on $L^2(D)$, so by 
the discussion above, there exists a symmetric Hunt process $X$ 
associated with $(\cE^{\vis},\cF^{\vis})$, taking values in $D$ 
with lifetime $\zeta$. We call $X$ a \emph{pure-jump process with visibility 
constraint} in $D$ associated with the 
kernel $k$. The process $X$ can be interpreted as the process obtained from the 
original pure-jump Markov process with jumping density $k$ by restricting its jumps from a point $x\in D$ to 
the visible area $D_x$ in $D$ from point $x$.

\subsection{Results}

Our first result is on comparability of the visibility constrained semi-norm to the semi-norm $|\cdot|_{H(D)}$, for a special class of kernels $k$ and domains $D$. From now on we will assume that, for some function $\ell:(0,\infty)\to(0,\infty)$ the following conditions are satisfied:
\begin{align}
&k(x,y)\asymp\frac{\ell(|y-x|)}{|y-x|^d}\qquad  (x,y\in D),\label{eq:K2intro}\\
&\int_0^{\infty}\left(r\wedge\frac 1 r\right)\ell(r)\d r<\infty,\label{eq:l1}\\
&\lambda^{-\gamma}\lesssim \frac{\ell(\lambda r)}{\ell(r)}\lesssim \lambda^d\qquad (\lambda\geq 1,\, r>0),\label{eq:l2}
\end{align}
for some constant $\gamma<2$. Note that the L\'{e}vy integrability condition \eqref{eq:l1} is equivalent to \eqref{eq:integrability1} and that \eqref{eq:l2} is a mild scaling condition ruling out fast decaying kernels. For the following comparability result we consider a special class of domains $D$, called uniform domains (see \autoref{def:uniform})
and jumping kernels such the function $\ell$ additionally satisfies a global scaling condition
\begin{align}
&\ell(\lambda r)\lesssim 
\lambda^{-\delta}\ell(r) \qquad (\lambda\geq 1,\, r > 
0)\label{eq:ell1} 
\end{align}
for some constant $0 < \delta \leq \gamma$. Given $s \in (0,1)$, the above conditions include examples like $\ell(r) = r^{-2s}$, $\ell(r) = r^{-2s} \ln(1+r^ {-1})^{\pm 1}$. 
   
\begin{theorem}\label{thm:comparable-intro}
Let $D$ be a bounded uniform domain. If the kernel $k$ is of the form \eqref{eq:K2intro} for a function $\ell$ satisfying \eqref{eq:l1}, \eqref{eq:l2} and \eqref{eq:ell1}, then 
\[
\cE^\cen(u,u)\lesssim \cE^\vis(u,u) \qquad (u\in L^2(D))\,.
\]
\end{theorem}

\begin{remark*}
The authors have been informed that, independently of this work, A.
Rutkowski has recently extended the comparability results \cite[(13)]{Dyd06} and \cite[Corollary 4.5]{PrSa17} allowing for a wider range of kernels and domains, see also \autoref{rem:comp-Rut}.
\end{remark*}

In \autoref{thm:comparability-weakly} we show a comparability result for a wider class of kernels including examples like $\ell(r) = \mathbbm{1}_{(0,1)}(r)$, under certain additional restrictions on the domain $D$.  \autoref{thm:comparable-intro} is a simple generalization of results from \cite{Dyd06} and \cite{PrSa17}, which cover bounded Lipschitz and bounded uniform domains and kernels $k$ comparable to the jumping density of the isotropic $2s$-stable L\'{e}vy process, $s\in(0,1)$, i.e.~case $\ell(r)=r^{-2s}$. The aforementioned papers show the comparability of a slightly weaker semi-norm to fractional Sobolev semi-norm $|\cdot|_{H^s(D)}$, see \eqref{eq:dyda}.

By applying \autoref{thm:comparable-intro}, we can recover useful density results and characterizations of spaces $\cF^\vis$ and $\widetilde{\cF}^\vis$. When $D$, $k$ and $\ell$ satisfy conditions in \autoref{thm:comparable-intro}, functions in $C_c(\overline{D})\cap\widetilde{\cF}^\vis$ are dense in $\widetilde{\cF}^\vis$, i.e.~$(\cE^\vis,\widetilde{\cF}^\vis)$ is a regular Dirichlet form on $L^2(\overline{D})$, see \autoref{rem:regularity}. Moreover, spaces $\cF^\vis$ and $\widetilde{\cF}^\vis$ are equal if and only if $\partial D$ is polar for the underlying isotropic unimodal L\'{e}vy process with radial jumping density $j(r)=\ell(r)r^{-d}$, \autoref{cor:prob-conseq}. A sufficient condition for $\cF^\vis = \widetilde{\cF}^\vis$ is provided in \autoref{cor:analysis-conseq}.

\medskip

Let us explain our main results regarding Poincar\'{e} inequalities. As explained above, we study these inequalities for dumbbell shaped domains on large scales. The main aim is to investigate the scaling behavior of the Poincar\'{e} constant with respect to large radii of balls. First of all, let us give a formal definition of the class of domains that we will study. 

\medskip

\hypertarget{ConditionA}{{\bf Condition A:}}
The set $D$ is of the form $D=D^+\cup\Gamma\cup D^-$, where $\Gamma$ is an open uniform set, $D^+$ and $D^-$ are disjoint open convex sets satisfying the following conditions: 
\begin{enumerate}[(i)]
	\item $\Gamma^* := \Gamma\setminus (D^+\cup D^-)$ is bounded,
	\item $|D^\pm\cap \Gamma|>0$, 
	\item There exists a $R_0>0$ such that for every $x_0 \in \Gamma^*$, $D^\pm \cap B(x_0,R) \asymp R^d$, for all $R\geq R_0$.
\end{enumerate}

\medskip

Sometimes we call the set $\Gamma$ the corridor. For $R\geq R_0$  and $x_0\in\Gamma^*$ set  
\begin{align}
\begin{split}\label{eq:setup2}
& D_R :=D\cap  B(x_0,R) ,\ D^\pm_R :=D^\pm\cap  B(x_0,R) , 
\end{split}
\end{align}
and $u_{D_R}=\frac{1}{|D_R|}\,\int_{D_R}u(x)\dx$. The following result is our main result concerning Poincar\'{e} inequalities, see also \autoref{thm:poincare-nonlocal}.

\begin{theorem}\label{thm:poincare-intro}
Assume $s \in (0,1)$ and $1 \leq p < d/s$. Let $D \subset \R^d$ be a domain 
satisfying \hyperlink{ConditionA}{Condition A}.
\begin{enumerate}[(i)]
	\item Then for all $R\geq R_0$ 
\begin{align*}
\int_{D_R}|u(x)-u_{D_R}|^p 
\dx \lesssim 
R^d\int_{D_R}\int_{D_{R, x}}\frac{|u(y)-u(x)|^p}{|x-y|^{d+sp}} \dy \dx 
\qquad (u\in  L^p(D_R)) \,.
\end{align*}
\item Let $s<1/p$. Assume that there exists a convex subset $\widetilde\Gamma$ of $\Gamma$ such that
\begin{align}\label{eq:Condition_C}
|\widetilde\Gamma\cap D_R^\pm|\gtrsim R \text{ for all }R\geq R_0.  
\end{align}
Then for all $R\geq R_0$
\begin{align*}
\int_{D_R}|u(x)-u_{D_R}|^p \dx\lesssim 
R^{d-1+sp}\int_{D_R}\int_{D_{R,x}}\frac{|u(y)-u(x)|^p}{|x-y|^{d+sp}}\dy \,
\dx \quad (u\in  L^p(D_R)) \,. 
\end{align*}
\end{enumerate}
\end{theorem}

\begin{figure}
	\centering
	\begin{subfigure}[b]{0.45\textwidth}
		\hspace*{4ex}\includegraphics[width=0.75\textwidth]{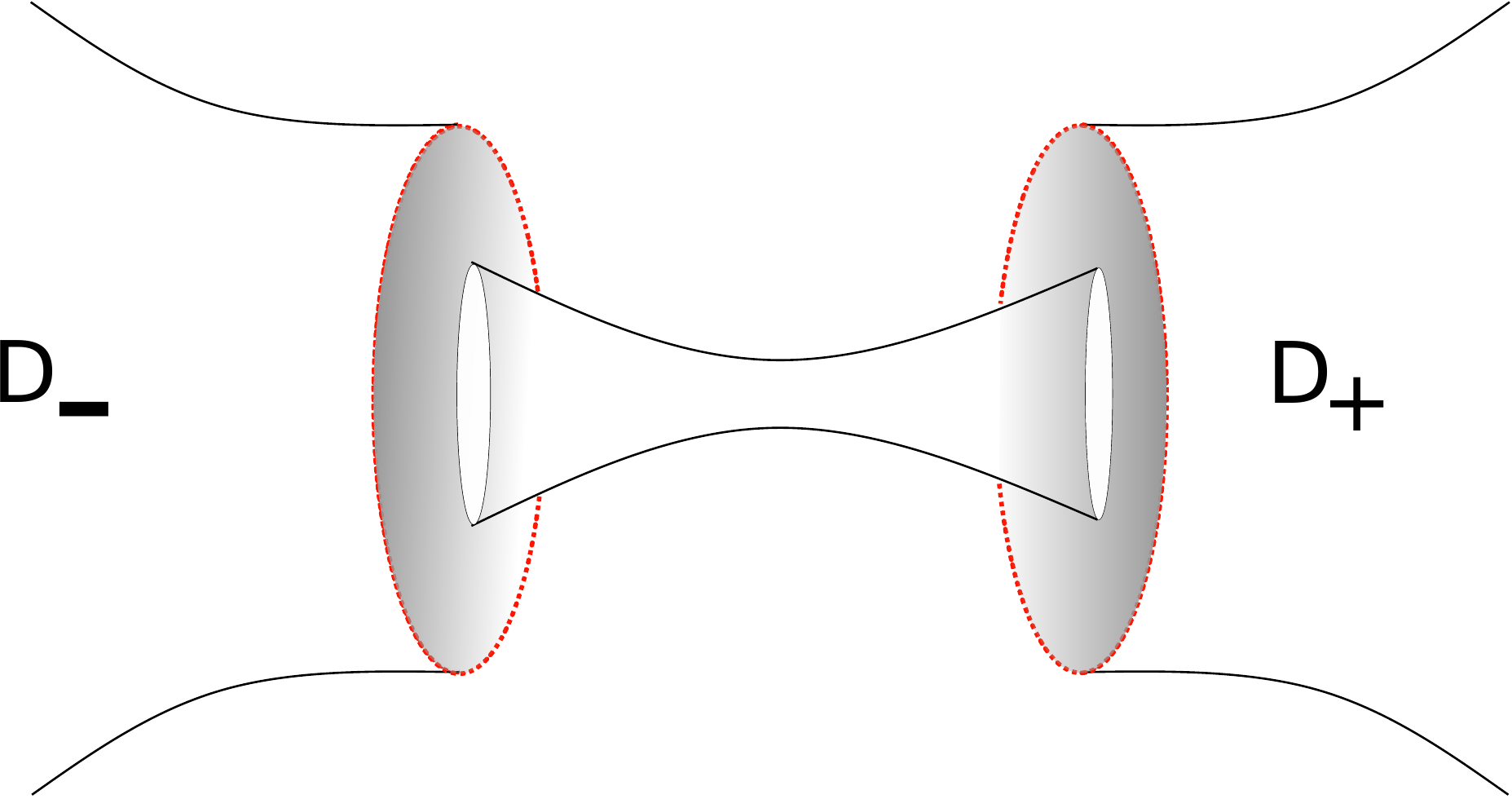}
		\caption{A domain satisfying \eqref{eq:Condition_C}}
		\label{fig:cond-A-straight}
	\end{subfigure}
	~ %add desired spacing between images, e. g. ~, \quad, \qquad, \hfill etc. 
	\begin{subfigure}[b]{0.45\textwidth}
		\hspace*{5ex}\includegraphics[width=0.7\textwidth]{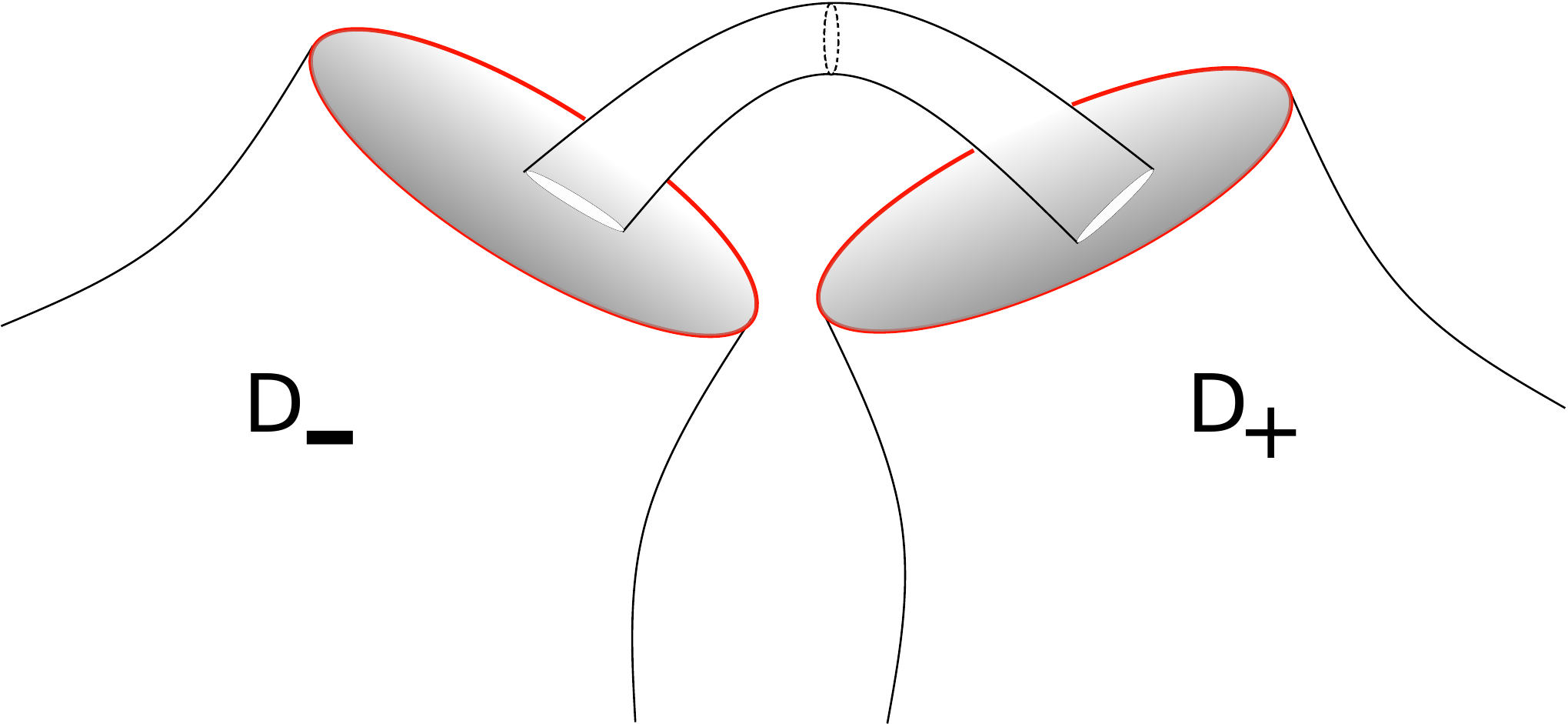}
		\caption{A domain violating \eqref{eq:Condition_C}}
		\label{fig:cond-A-curved}
	\end{subfigure}
	\caption{Two unbounded domains satisfying Condition A}\label{fig:exampes-cond-A}
\end{figure}

Part (ii) describes a remarkable phenomenon. If the corridor $\Gamma$ allows large regions of $D^-$ and $D^+$ to be connected by long jumps, then this geometric situation allows for smaller Poincar\'{e} constants in the case of small values of $s$. We provide the proof of \autoref{thm:poincare-intro} in \autoref{sec:poincare} together with a discussion in which sense the result is sharp.     

\subsection{Organization of the article and notation}

\emph{Organization:} The article is organized as follows. In \autoref{sec:comparability} we prove \autoref{thm:comparable-intro} and discuss the consequences of this result for the corresponding Markov process. In \autoref{ex:counter-weak-kernels} we provide a counterexample of a weakly singular kernel for which \autoref{thm:comparable-intro} does not hold. For a more restricted class of domains we are able to formulate and prove a version of \autoref{thm:comparable-intro} that allows for weakly singular kernels, cf. \autoref{thm:comparability-weakly}. \autoref{sec:poincare} is devoted to the proof and the discussion of the Poincar\'{e} inequality in different settings, in particular of \autoref{thm:poincare-intro}. 

\medskip

\emph{Notation:} Throughout the paper, we use the notation $f\lesssim g$ ($f\gtrsim g$) if there exists a constant $c>0$ such that $f(x)\leq c g(x)$ ($f(x)\geq c g(x)$) for every $x$. We write $f\asymp g$ if $f\lesssim g$ and $f\gtrsim g$. 

\medskip

\emph{Acknowledgement:} The authors would like to thank Krzysztof Bogdan, Bart{\l}omiej Dyda and Armin Schikorra for helpful discussions on the subject of this work.

\section{Comparability of bilinear forms}\label{sec:comparability}

The first aim of this section is to prove \autoref{thm:comparable-intro}. In order to do so, we first discuss the notion of uniform domains. After the proof of \autoref{thm:comparable-intro} we comment on a more refined result that would follow with similar techniques. In \autoref{ex:counter-weak-kernels} we provide an example showing that conditions \eqref{eq:K2intro}, \eqref{eq:l1}, and \eqref{eq:l2} are not sufficient for \autoref{thm:comparable-intro}. Next, we discuss consequences of \autoref{thm:comparable-intro} regarding the Potential Theory related to visibility constrained jump processes. Finally, we provide a comparability result, \autoref{thm:comparability-weakly}, that allows for weakly singular kernels.

\medskip
	
Let $\ell:(0,\infty)\to(0,\infty)$ be a function satisfying \eqref{eq:l1}, $j(r):=\frac{\ell(r)}{r^d}$, $r>0$, and
\begin{align*}
 \psi(|\xi|):=\int_{\R^d}(1-\cos(\xi\cdot x))j(|x|)\dx,\qquad \xi\in\R^d.
\end{align*}
Since \eqref{eq:l1} is equivalent to
\begin{align}\label{eq:j}
\int_0^\infty (1\wedge r^2)j(r)r^{d-1}\d r<\infty,
\end{align}
the function $\psi$ is the characteristic exponent of an isotropic pure-jump L\'{e}vy process with a radial L\'{e}vy density $j$. If, additionally, $\ell$ is a non-increasing function, the corresponding process is an isotropic unimodal L\'{e}vy process. Our first comparability result concerns jumping kernels $k$ which are comparable to jumping kernels of isotropic unimodal L\'{e}vy processes satisfying a scaling condition \eqref{eq:ell1} and bounded uniform domains. 

\begin{definition}
Let $D$ be a domain and $\mathcal W$ a Whitney decomposition of $D$, see \cite[Chapter I.2.3]{JoWa84} for details. For $Q,S\in\mathcal W$ denote by $D(Q,S)$ the long distance between cubes $Q$ and $S$, defined by $D(Q,S):=\LL(Q)+\dist (Q,S)+\LL(S)$, where $\LL(Q)$ is the side length of cube $Q$. A chain $[Q,S]$ of size $k\in\mathbb N$ connecting cubes $Q, S \in\mathcal W$ is a series of 
cubes $\{Q_1,...,Q_k\}$ in $\mathcal W$ such that $Q_1=Q$, $Q_k=S$ and $Q_i$ 
and $Q_{i+1}$ touch each other for all $i$.  Let $\varepsilon>0$. A chain $[Q, S]$ 
is $\varepsilon$-admissible if
\begin{enumerate}[(i)]
\item the length of the chain is bounded by
\[
\LL([Q,S]):=\sum_{i=1}^k\LL(Q_i)\leq\frac 1\varepsilon D(Q,S)
\]
\item there exists $j_0\leq k$ such that the cubes in the chain satisfy
\[
\LL(Q_j)\geq\begin{cases}
\varepsilon D(Q,Q_j) & j\leq j_0\\
\varepsilon D(Q_j,S) & j\geq j_0.\\
\end{cases}
\]
For an admissible chain $[Q,S]$ we denote the central cube $Q_{j_0}$ as $Q_S$.
\end{enumerate}
\end{definition}

Note that by choosing cubes of smaller size in the Whitney decomposition one can get that $S\subset\cup_{x\in Q} D_x$ for all cubes $Q,S\in\mathcal W$ that touch each other. 

\begin{definition}\label{def:uniform}
We say that a domain $D\subset\R^d$ is a uniform domain if there exists a $\varepsilon>0$ and a Whitney covering $\mathcal W$ of $D$ such that for any pair of cubes $Q, S \in\mathcal W$, there exists an $\varepsilon$-admissible chain $[Q, S]$. 
\end{definition}

\begin{remark*}\label{rem:domains} 
Let $k$ be the kernel satisfying \eqref{eq:K2intro} for $\ell(r)=r^{-2s}$, $0<s<1$, i.e. a kernel comparable to the L\'{e}vy density $j(r)=r^{-d-2s}$ of the isotropic $2s$-stable L\'{e}vy process. If one of the following cases holds
\begin{itemize}
 \item $D$ is the domain above the graph of a Lipschitz function;
 \item $D$ is a connected component of the complement of a bounded 
Lipschitz open set;
 \item $D$ is a bounded uniform domain,
\end{itemize}
then by \cite[(13)]{Dyd06} and \cite[Corollary 4.5]{PrSa17}  
there exists a constant $c_1=c_1(d,D,s)>0$ such that for all 
$u:D\to\R$

\begin{align}
 \cE^\cen(u,u)&\leq c_1 
\int\limits_D\int\limits_{B(x,\delta_D(x)/2)}(u(y)-u(x))^2k(x,y)\dy \,
\dx\label{eq:dyda} \\ 
 &\leq c_1\cE^\vis(u,u)\nonumber,
\end{align}
where $\delta_D(x)=\text{dist}(x,\partial D)$. This proves 
$\cE^\vis\asymp \cE^\cen$ on $L^2(D)\times L^2(D)$, $\cF^\vis=\cF^\cen$ and $\widetilde{\cF}^\vis=\widetilde{\cF}^\cen$.
\end{remark*}

\autoref{thm:comparable-intro} is the extension of this comparison result for a wider class of kernels, satisfying \eqref{eq:K2intro}, \eqref{eq:l1}, \eqref{eq:l2} and \eqref{eq:ell1}. In the proof we follow the approach in \cite{PrSa17}. 

\proofof \autoref{thm:comparable-intro}:
Throughout the proof we use the semi-norm in the duality form
\[
 |f|_{H(D)}\asymp\sup_{||g||_{
L^2(D \times D)}\leq 1} \int_D \int_D 
|f(x)-f(y)| \frac{\ell(|x-y|)^{1/2}}{|x-y|^{d/2}}g(x,y)\dy \dx.
\]

By construction, there exists a constant $c_1>0$ such that $Q^*:=(1+c_1)Q\subset D$ for all $Q \in \cW$. Using the Whitney decomposition, we divide the semi norm into two parts, 
\begin{align*}
 &\int_D \int_D 
|f(x)-f(y)|\frac{\ell(|x-y|)^{1/2}}{|x-y|^{d/2}}g(x,y)\dy \d 
x\\
&=\sum_{Q\in\cW}
\int_Q \int_{Q^*} 
|f(x)-f(y)|\frac{\ell(|x-y|)^{1/2}}{|x-y|^{d/2}}g(x,y)\dy \dx\\  
 &\quad +\sum_{Q,S\in\cW}\int_Q \int_{S\setminus Q^*} 
|f(x)-f(y)| \frac{\ell(|x-y|)^{1/2}}{|x-y|^{d/2}}g(x,y)\dy \dx=: 
J_1+J_2.  
\end{align*} 
Since $Q^*\subset D_x$ for all $x\in Q$, by H\"{o}lder's inequality we immediately deduce
\[
 J_1 \leq \left(\int_D \int_{D_x} 
(f(x)-f(y))^2\frac{\ell(|x-y|)}{|x-y|^d}\dy \dx\right)^{1/2}.
\]
For the next term, note that $|x-y|\asymp D(Q,S)$ for $x\in Q$ and $y\in S\setminus Q^*$. For a cube $P$ in an admissible chain, we denote by $\mathcal N(P)$ the following cube in the same chain. Applying the triangle inequality along the chain $[Q,Q_S)$ and taking into account \eqref{eq:l2}, we obtain
\begin{align*}
 J_2&\lesssim\sum_{Q,S\in\cW}\int_Q \int_{S}  
|f(x)-f_Q|\frac{\ell(D(Q,S))^{1/2}}{D(Q,S)^{d/2}}g(x,y)\dy \dx\\
 & \quad +\sum_{Q,S\in\cW}\int_Q \int_{S} 
\sum_{P\in[Q,Q_S)}|f_P-f_{\mathcal 
N(P)}|\frac{\ell(D(Q,S))^{1/2}}{D(Q,S)^{d/2}}g(x,y)\dy \dx\\ 
&\quad +\sum_{Q,S\in\cW}\int_S \int_{Q} 
\sum_{P\in[S,Q_S)}|f_P-f_{\mathcal 
N(P)}|\frac{\ell(D(Q,S))^{1/2}}{D(Q,S)^{d/2}}g(x,y)\dx \dy\\ 
 &\quad +\sum_{Q,S\in\cW}\int_S \int_{Q} 
|f_{S}-f(y)|\frac{\ell(D(Q,S))^{1/2}}{D(Q,S)^{d/2}}g(x,y)\dx \,
\dy=:I_1+I_2+I_3+I_4. 
\end{align*} 
In the following calculations we will frequently apply the following essential property of the Whitney decomposition, from \cite[Lemma 3.13]{PrTo15}: for all $b>a>d-1$
\begin{align}\label{eq:whitney1}
 \sum_{S\in\mathcal W}\frac{\LL(S)^a}{D(Q,S)^b}\lesssim \LL(Q)^{a-b}\qquad (Q\in\mathcal W).
\end{align}
By H\"{o}lder's inequality, we deduce
\begin{align*}
 I_1^2&\lesssim \sum_{Q,S\in\cW}\int_Q \int_{S} 
(f(x)-f_Q)^2\frac{\ell(D(Q,S))}{D(Q,S)^{d}}\dy \dx\\ 
 &\leq \sum_{Q,S\in\cW}\int_Q \int_{S}\frac{1}{\LL(Q)^d}\left(\int_Q 
(f(x)-f(z))^2\dz\right)\frac{\ell(D(Q,S))}{D(Q,S)^{d}}\dy \dx\\ 
 &= \sum_{Q\in\cW}\left(\int_Q \int_Q \frac{(f(x)-f(z))^2}{\LL(Q)^d} \dz\dx 
\cdot\sum_{S\in\cW}\frac{\LL(S)^d\ell(D(Q,S))}{D(Q,S)^{d}}\right)\\ 
&\overset{\eqref{eq:ell1}}{\lesssim} \sum_{Q\in\cW}\left(\int_Q \int_Q 
\frac{(f(x)-f(z))^2}{\LL(Q)^d}\ell(\LL(Q)) dz\dx \cdot 
\LL(Q)^{\delta}\sum_{S\in\cW}\frac{\LL(S)^d}{D(Q,S)^{d+\delta}}
\right)\\ 
&\underset{\eqref{eq:ell1}}{\overset{\eqref{eq:whitney1}}{\lesssim}}\sum_{Q\in\cW}\int_Q \int_Q 
\frac{(f(x)-f(z))^2}{|z-x|^d}\ell(|z-x|) dz\dx\\
&\lesssim\int_D  
\int_{D_x} \frac{(f(x)-f(z))^2}{|z-x|^d}\ell(|z-x|) dz\dx.  
\end{align*}
 Note that, by the properties of the Whitney covering, there exists a constant $c_2>0$ depending only on the covering $\mathcal W$ such that $\mathcal N(P)\subset \tilde P\cap D_{P}$, where $D_P:=\cap_{x\in P}D_x$ and $\tilde P:=(1+c_2)P$. Recall also that $\LL(P)\asymp \LL(\mathcal N(P))$, since the cubes $P$ and $\mathcal N(P)$ touch. Furthermore, we note that there exists  $\rho>0$, such that for all $Q,S\in\mathcal W$ and all $P\in[Q,Q_S]$, $Q\subset B\left(x_P,\rho\, \LL(P)\right)$ and $D(Q,S)\asymp D(P,S)$. Here $x_P$ is the center of cube $P$. Also, we write $Q\leq P$ if there exists $S\in\mathcal W$ such that $P\in[Q,S]$. Therefore, 

\begin{align*}
 I_2&\lesssim \sum_{Q,S\in\cW}\sum_{P\in[Q,Q_S)}\int_Q \int_{S} 
\int_P\int_{\mathcal N(P)}\frac{|f(z)-f(w)|}{\LL(P)^{2d}}\frac{\ell(D(Q,S))^{1/2}}{D(Q,S)^{d/2}}g(x,y)\dw\dz\dy \dx\\ 
&\leq \sum_{P\in\mathcal W}\int_P\int_{\tilde P\cap D_P}\frac{|f(z)-f(w)|}{\LL(P)^{2d}}\dw\dz\sum_{Q\leq P}\sum_{S\in \cW}\int_Q \int_{S} 
\frac{\ell(D(Q,S))^{1/2}}{D(Q,S)^{d/2}}g(x,y)\dy \dx\\ 
&\leq \sum_{P\in\mathcal W}\LL(P)^{-2d}\int_P\int_{\tilde P\cap D_P}|f(z)-f(w)|\dw\dz\cdot \sum_{Q\leq P}\int_Q \left(\int_D g(x,y)^2\dy\right)^{1/2}\dx \\
&\quad \cdot\left(\sum_{S\in \cW} 
\frac{\LL(S)^d\ell(D(P,S))}{D(P,S)^{d}}\right)^{1/2}\\ 
&\overset{\eqref{eq:ell1}}{\lesssim} \sum_{P\in\mathcal W}\LL(P)^{-2d}\int_P\int_{\tilde P\cap D_P}|f(z)-f(w)|\dw\dz\cdot \sum_{Q\leq P}\int_Q G(x)\dx \\
&\quad \cdot\left(\ell(\LL(P))\LL(P)^{\delta}\sum_{S\in \cW} 
\frac{\LL(S)^d}{D(P,S)^{d+\delta}}\right)^{1/2}\\
&\overset{\eqref{eq:whitney1}}{\lesssim} \sum_{P\in\mathcal W}\LL(P)^{-2d}\ell(\LL(P))^{1/2}\int_P\int_{\tilde P\cap D_P}|f(z)-f(w)|\dw\dz\cdot \sum_{Q\leq P}\int_Q G(x)\dx, 
\end{align*}
where $G(x):=\left(\int_D g(x,y)^2\dy\right)^{1/2}$. By H\"{o}lder's inequality and \cite[Lemma 3.11]{PrTo15} one obtains
\begin{align*}
 \sum_{Q\leq P}\int_Q G(x)\dx\leq  \sum_{Q\leq P}\LL(Q)^{d/2}\left(\int_D \int_D g(x,y)^2\dy\dx\right)^{1/2}\lesssim \LL(P)^{d/2}.
\end{align*}
Therefore, by applying H\"{o}lder's inequality once more, we obtain
\begin{align*}
 I_2&\lesssim \sum_{P\in\mathcal W}\LL(P)^{-3d/2}\ell(\LL(P))^{1/2}\int_P\int_{\tilde P\cap D_P}|f(z)-f(w)|\dw\dz\\
 &\leq \sum_{P\in\mathcal W}\LL(P)^{-d/2}\ell(\LL(P))^{1/2}\left(\int_P\int_{\tilde P\cap D_P}(f(z)-f(w))^2\dw\dz\right)^{1/2}\\
&\overset{\eqref{eq:ell1}}{\lesssim} \sum_{P\in\mathcal W} \left(\int_P\int_{\tilde P\cap D_P}(f(z)-f(w))^2\frac{\ell(|z-w|)}{|z-w|^{d}}\dw\dz\right)^{1/2}\\
&\leq \left(\int_D\int_{D_z}(f(z)-f(w))^2\frac{\ell(|z-w|)}{|z-w|^{d}}\dw\dz\right)^{1/2}.\\
\end{align*}
By applying analogous calculations to terms $I_3$, $I_4$ and combining all established estimates, the proof is concluded. 
\qed

\begin{remark}\label{rem:comp-Rut}
Note that by following the proof of \cite[Lemma 4.1., Lemma 4.3]{PrSa17} one can obtain a stronger result of the form \eqref{eq:dyda}, i.e.
\begin{align}
 \int\limits_D\Bigg(\int\limits_{D}|u(y)-u(x)|^q k(x,y)\dy\Bigg)^{p/q} 
\!\!\!\!\dx&\lesssim 
\int\limits_D\Bigg(\int\limits_{B(x,\delta_D(x)/2)}\!\!\!|u(y)-u(x)|^q k(x,y)\dy\Bigg)^{p/q} \!\!\!\!\!\!\!\dx\nonumber,
\end{align}
for $p,q>1$ and the kernel $k$ of the form 
\begin{align*}
 k(x,y)\asymp\frac{\ell(|x-y|)}{|x-y|^d}\qquad (x,y\in D),
\end{align*}
where the function $\ell:(0,\infty)\to(0,\infty)$ is non-increasing and satisfies
\begin{align*}
\ell(\lambda r)\lesssim \lambda^{-qs}\ell(r)\qquad (\lambda\geq1,\, r>0)
\end{align*}
for some $d(1/p-1/q)_+<s<1$. This approach has recently been pursued by A. Rutkowski.
\end{remark}

The following example shows that one cannot expect \eqref{eq:dyda} resp.~the result of the aforementioned remark to hold for general kernels $k$ satisfying only \eqref{eq:K2intro}, \eqref{eq:l1} and \eqref{eq:l2}, no matter how regular the domain $D$ is. 

\begin{example}\label{ex:counter-weak-kernels}
	Let $D=[0,1]^2$, $k(x,y)=\frac{1}{|y-x|^2}$. Define a sequence $(u_n)$ in $L^2(D)$ by 
	\begin{align*}
	u_n(x_1, x_2)=\mathbbm{1}_{(0,1/n)}(x_1+x_2), \text{ where } n\in\N.
	\end{align*}
    Then \eqref{eq:dyda} fails because, as we will show,
    \begin{align}\label{eq:show}
\frac{ \int_D\int_{B(x,\sfrac{\delta_D(x)}{2})}(u_n(y)\!-\!u_n(x))^2k(x,y)\dy \, \dx}{\cE^\cen (u_n, u_n) } \longrightarrow 0 \quad \text{ as } n \to \infty\,.
    \end{align}
Since $D$ is convex, the two quantities $\cE^\cen$ and $\cE^\vis$ are equal. Note
	\begin{align*} 
	&\int\limits_D\int\limits_{B(x,\sfrac{\delta_D(x)}{2})}(u_n(y)-u_n(x))^2k(x,y)\dy \,
	\dx \leq \int_{A_n}\int_{(D\setminus A_n)\cap B(0,\frac{1}{4n})}\frac{1}{|y-x|^2}\dy \dx,
	\end{align*}
	where $A_n=\{(x_1,x_2)\in(0,1)^2:x_1+x_2<1/n\}$. Furthermore, 
	\begin{align*}
	&\int_{A_n}\int_{(D\setminus A_n)\cap B(0,\frac{1}{4n})}\frac{1}{|y-x|^2}\dy \dx\lesssim 
	\int_{A_n}\int_{\dist(x,D\setminus A_n)}^{\frac{1}{4n}}\frac 1 s \,\d s\, \dx\\
	&= -\frac{\log{4}}{2n^2} -\frac{\log{n}}{2n^2}-\int_{A_n}  \log(\dist(x,D\setminus A_n))\,\dx\\
	&\leq -\frac{\log{n}}{2n^2}-\int_0^{1/n}\int_0^{1/n-x_1}  \log\left(\frac 1 n -x_1 - x_2\right)\,\dx_2\,\dx_1\\
	&=-\frac{\log{n}}{2n^2}-\left(-\frac{\log{n}}{2n^2}-\frac 3{4n^2}\right)\leq \frac 1{n^2}\,. 
	\end{align*} 
On the other hand, 
	\begin{align*}
	\cE^\cen(u_n,u_n) = \int_{A_n}\int_{D\setminus A_n}\frac{1}{|y-x|^2}\dy \dx\gtrsim 
	\int_{A_n}\int_{\frac{\sqrt{2}}2\dist(x,D\setminus A_n)}^{1}\frac 1 s \,\d s\, \dx\gtrsim 
	\frac 
	{\log{n}}{n^2}\,,
	\end{align*}
which together with the previous estimate implies \eqref{eq:show}. 
\end{example}

\medskip

By \cite[Corollary 23]{BGR14}, conditions of \autoref{thm:comparable-intro} imply the following global estimates for the L\'{e}vy density $j$ in terms of the radial non-decreasing majorant $\psi^*$ of the characteristic exponent $\psi$,
\begin{align*}
j(r)\asymp \frac{\psi^*(r^{-1})}{r^d},\qquad (r>0),     
\end{align*}
i.e. $\ell(r)\asymp \psi^*(r^{-1}):=\sup
_{0\leq u\leq r^{-1}}\psi(u)$. Combining \autoref{thm:comparable-intro} with the results on boundary behavior of the censored process, see \cite{BBC03} for the stable case 
and \cite{Wag18} for the general case, we arrive to the following corollary.  

\begin{corollary}\label{cor:prob-conseq} Let $D$ be an open bounded uniform 
domain and kernel $k$ such that \eqref{eq:K2intro} holds for a function $\ell:(0,\infty)\to(0,\infty)$ satisfying \eqref{eq:l1}, \eqref{eq:l2} and \eqref{eq:ell1}. The following 
statements are equivalent
 \begin{enumerate}[(i)]
  \item $X$ is recurrent and therefore conservative, $\P_x(X_{\zeta-}\in\partial D,\zeta<\infty)=0$;
  \item $\partial D$ is polar for the L\'{e}vy process with the characteristic exponent $\psi$;
  \item $1\in \cF^\vis$;
  \item $1\in \cF^\cen$;
    \item $\cF^\vis=\widetilde{\cF}^\vis$.
 \end{enumerate}
\end{corollary} 

As a direct consequence of \autoref{cor:prob-conseq} (see also \cite{BBC03}, \cite{Wag18}) we get sufficient conditions (in terms of $\delta$ and $\gamma$) on the equivalence of spaces $\cF^\vis$ and $\widetilde{\cF}^\vis$ when $D$ is a bounded Lipschitz domain. Note that the boundary $\partial D$ is polar for the underlying unimodal L\'{e}vy process if and only if it is of zero capacity, i.e.
\begin{align*}
 \text{Cap}_\psi(\partial D):=\inf\limits_{f,U}\left\{\cE_1^{\R^d}(f,f) |\, f\in L^2(\R^d), \, f\geq 1 \text{ a.e. on }U,\, \partial D\subset U \text{open}\right\}=0,
\end{align*}
where $\cE^{\R^d}(f,f):=\int_{\R^d}\int_{\R^d}(f(x)-f(y))^2\frac{\ell(|x-y|)}{|x-y|^d}\dy\dx$, see for example \cite[Section II.3.]{Ber96} and \cite{FOT94}. Under conditions $\eqref{eq:l2}$ and $\eqref{eq:ell1}$ we get the lower and upper bound on $\text{Cap}_\psi$ in terms of the Riesz capacity of order $n-\delta$ and $n-\gamma$ respectively, 
\begin{align*}
 \text{Cap}_{n-\delta}(\partial D)\lesssim \text{Cap}_\psi(\partial D)\lesssim \text{Cap}_{n-\gamma}(\partial D).
\end{align*}
By using the well known relation between the Riesz capacity and the Hausdorff dimension of a set, see e.g.~\cite{adams}, we arrive to the following result.
\begin{corollary}\label{cor:analysis-conseq} 
	Let $D$ be an open bounded Lipschitz domain and kernel $k$ such that \eqref{eq:K2intro} holds for a function $\ell:(0,\infty)\to(0,\infty)$ satisfying \eqref{eq:l1}, \eqref{eq:l2} and \eqref{eq:ell1}. Then
 \begin{enumerate}[(i)]
  \item if $\gamma\leq 1$ then $\cF^\vis=\widetilde{\cF}^\vis$,
  \item if $\delta> 1$ then $\cF^\vis\subsetneq\widetilde{\cF}^\vis$.
 \end{enumerate}
\end{corollary}
For a more general version of this result, stated in terms of the Hausdorff dimension of the boundary of an open uniform domains and for spaces $\cF^\cen$, $\widetilde{\cF}^\cen$, see \cite[Corollary 1.3]{Wag18} and \cite[Corollary 2.8]{BBC03}.

\begin{remark}\label{rem:regularity}
$\quad$
\begin{enumerate}[(i)]
\item Note that the \cite[Corollary 2.6]{BBC03} and \cite[Corollary 2.9]{Wag18}, which we apply here, are both stated for bounded open $d$-sets $D$ in $\R^d$, i.e. for sets such that
\begin{align*}
 |B(x,r)\cap D|\gtrsim r^d \qquad (x\in \overline{D},\,0<r<1).
\end{align*}
By \cite[II.1.1. Example 4]{JoWa84} and \cite[p.2496]{PrSa17}, every uniform domain in $\R^d$ satisfies this condition.
\item Let $D$ be a bounded Lipschitz set and $\psi(r) = 
r^{2s}$ for $s \in 
(0,1)$. As a consequence of \autoref{cor:prob-conseq}, the visibility constrained process $X$ is recurrent if 
$s\in(0,\frac12]$, otherwise it is transient. This also follows directly from \cite[Corollary 10, Corollary 11]{Dyd06} and \cite[Theorem 1.1]{BBC03}.
\item Another consequence of \autoref{cor:prob-conseq} is that for open bounded uniform domains $D$ the Dirichlet form $(\cE^\vis, \widetilde{\cF}^{\vis})$ is regular on $L^2(D)$ with the core $C_c(\overline D)\cap\widetilde{\cF}^\vis$.
\end{enumerate}
\end{remark}

The scaling condition \eqref{eq:ell1}, which allowed for the application of inequality \eqref{eq:whitney1}, was important in \autoref{thm:comparable-intro} for treating admissible paths of arbitrary sizes, which are characteristic for uniform domains. By posing additional constraints on the size of the admissible paths, we can extend this result to a wider class of kernels $k$, allowing for weakly singular kernels. To this end, we introduce the following condition on domains $D\subset\R^d$: 

\medskip

\hypertarget{ConditionB}{{\bf Condition B:}} There exists a constant 
$N=N(D)\in \N$ and $c=c(D)\geq1$ such that for almost every $x\in D$ 
and almost every $y\not\in D_x$ there exists a $k\leq N$ and cubes $Q_1,...,Q_k$ in 
$D$ such that
\begin{itemize}
 \item $\LL(Q_i)\asymp |x-y|$, for all $i=1,...,k$,
 \item $Q_1\subset D_x$, $Q_k\subset D_y$ and $\dist(x,B_1),\dist(y,B_k)\asymp|x-y|$, 
 \item $Q_{i+1}\subset D_{Q_i}=\cap_{x\in Q_i} D_x$, $\dist(Q_i,Q_{i+1})\asymp|x-y|$, 
\end{itemize}
where the constants in the comparisons depend only on $D$. We call this family of cubes an admissible path of length $k$ for $x\in D$ and $y\in D_x^c$.

\begin{remark}

\begin{enumerate}[(i)]
 \item Note that this condition is satisfied when $D$ is a uniform domain with admissible chains of bounded size. Furthermore, a domain $D$ satisfying \hyperlink{ConditionB}{Condition B} is a uniform domain (see \cite{PrSa17} for the discussion on equivalent definitions of uniform domains). 
  \item One can easily show that a connected finite union of open bounded convex sets $K_i$ satisfies \hyperlink{ConditionB}{Condition B} if for every two components $K_i$ and $K_j$
  \begin{align*}
   \overline{K_i}\cap\overline{K_j}=\emptyset\qquad \text{or}\qquad |K_i\cap K_j|>0.
  \end{align*}
  \item An example of a uniform domain that does not satisfy \hyperlink{ConditionB}{Condition B} is a Koch snowflake domain, see \cite{Jon81}.
 \end{enumerate}
\end{remark}

\begin{theorem}\label{thm:comparability-weakly}
Let $D$ be an open set in $\R^d$ satisfying \hyperlink{ConditionB}{Condition B} and $k(x,y)\asymp j(|y-x|)$, $x,y\in D$, for some function $j:(0,\infty)\to(0,\infty)$ satisfying \eqref{eq:j} and
\begin{align}\label{eq:g}
g_1(\lambda)j(r)\lesssim j(\lambda r)\lesssim g_2(\lambda)j(r),\qquad (\lambda\geq 1,\, r>0)
\end{align} 
for some non-decreasing functions $g_i:[1,\infty)\to (0,\infty)$, $i=1,2$. Then 
\begin{align*}
 \cE^\cen(u,u)\lesssim \cE^\vis(u,u) \qquad (u\in L^2(D)) \,.
\end{align*}
\end{theorem}

\proof 
\begin{align*}
 &\cE^\cen(u,u)-\cE^\vis(u,u)=\int_D\int_{D_x^c} (u(y)-u(x))^2 k(x,y)\dy \, 
\dx\nonumber\\ 
 &\lesssim 
\int_D\int_{D_{x}^c}\frac{1}{|x-y|^{kd}}\int_{Q^{x,y}_1}\int_{Q^{x,y}_2}
...\int_{Q^{x,y}_{k}}(u(x)-u(y))^2 j(|x-y|)\dx_{k}...\d 
x_{1}\dy  
\dx\nonumber\\ 
&\overset{\eqref{eq:g}}{\hspace*{-2ex}\lesssim} 2^N\int_D\int_{D_{x_0}^c}\sum_{i=0}^k\int_{Q^{x,y}_1}
...\int_{Q^{x,y}_{k}}\!\!\! \frac{(u(x_{i+1})-u(x_{i}))^2}{|x-y|^{kd}} 
j(|x_{i+1}-x_{i}|)\dx_{k}...\d 
x_{0}
\nonumber\\  
 &\lesssim \int_D\int_{D_{x}^c}\sum_{i=2}^{k} 
\int_{Q^{x,y}_i}\int_{Q^{x,y}_{i-1}} (u(x_i)-u(x_{i-1}))^2 
\frac{j(|x_{i} -x_{i-1}|)}{|x_{i}-x_{i-1}|^{2d}}\dx_{i-1} \dx_{i}\dy \d 
x\nonumber\\
 &\quad + \int_D\int_{D_{x}^c}\int_{Q^{x,y}_1} (u(x_1)-u(x))^2 
\frac{j(|x_{1}-x|)}{|x_{1}-x|^{d}} \dx_{1}\dy \dx \\ 
 &\quad +\int_D\int_{D_{x}^c}\int_{Q^{x,y}_{k}} 
(u(y)-u(x_{k}))^2 \frac{j(|y-x_{k}|)}{|y-x_{k}|^{d}} \dx_{k}\dy \d 
x\\ 
&\lesssim  \int_D\int_{D_{x}^c}\sum_{i=2}^{k}  
\iint_{A^{x,y}_i}\frac{ (u(z)-u(w))^2 j(|z-w|)}{|z-w|^{2d}} \d w\, \d z\, \d 
y\,  
\dx \\
&\quad + \int_D\int_{D_{x}^c}\int_{B^{x,y}} (u(z)-u(x))^2  
\frac{j(|z-x|)}{|z-x|^{d}} \d z\dy \dx \\
 &\quad +\int_D\int_{D_{x}^c}\int_{C^{x,y}} (u(y)-u(z))^2 
\frac{j(|y-z|)}{|y-z|^{d}} \d z\dy \dx \\
&=I_1+I_2+I_3
\end{align*}
where 
\begin{align*}
 A^{x,y}_i&:=\{(z,w)\in Q^{x,y}_i\times D_z \big| \,c_1^{-1}|z-w| \leq |v-z|\leq c_1|z-w| \text{ for }v=x,y\}\\
 A^{x,y}_i&\supset Q^{x,y}_i\times Q^{x,y}_{i-1},\\ 
 B^{x,y}&:=\{z\in D_x \big|\, c_2^{-1}|z-x|\leq|y-x|\leq 
c_2|z-x|\}\supset Q^{x,y}_1 \\
 C^{x,y}&:=\{z\in D_y \big|\ c_2^{-1}|z-y|\leq|y-x|\leq c_2|z-y|\}\supset Q^{x,y}_{k}, 
 \end{align*}
for some constants $c_1,c_2\geq 1$ depending only on $D$. Since 
$A_i^{x,y}$ are mutually  disjoint and $\cup_i A_i^{x,y}\subset A^{x,y}:=\{(z,w)\in D\times D_z|c_1^{-1}|z-w| \leq |v-z|\leq c_1|z-w|,\ \text{for }v=x,y\}$ it 
follows that the integral $I_1$ is less than
\begin{align*}
&c \int_D\int_{D_{x}^c}\iint_{A^{x,y}}\frac{ (u(z)-u(w))^2 
j(|z-w|)}{|z-w|^{2d}} \d w\, \d z\, \dy \, \dx,\\ 
&\leq 
\int_D\int_{D_z}\frac{ 
(u(z)-u(w))^2 j(|z-w|)}{|z-w|^{2d}}\int_{B(z,c_1|w-z|)}\int_{B(z,c_1|w-z|)}\dy \dx \dz \dw\\
&\asymp \cE^\vis(u,u)  
\end{align*}
Similarly,
\begin{align*}
I_2\leq 
\int_D\int_{D_x}\int_{B(x,c_2|z-x|)}(u(z)-u(x))^2 
\frac{j(|z-x|)}{|z-x|^{d}}\dy \d z\dx\asymp\cE^\vis(u,u) \,.
\end{align*}
The same inequality follows for the integral $I_3$.
\qed

\vspace{1cm}

\section{Poincar\'{e} inequality of the visibility constrained bilinear form}\label{sec:poincare}

The aim of this section is to prove \autoref{thm:poincare-intro} and related results. An interesting example of a domain where the classical Poincar\'{e} inequality fails is given by a dumbbell shaped manifold, e.g., see \cite[Example 2.1]{Sal09}. Such domains can be decomposed into a disjoint 
union $D^-\cup \Gamma\cup D^+$, where $D^-$, $D^+$ each are isometric to the outside of some compact domain with smooth boundary in $\R^d$ and $\Gamma$ is a smooth compact manifold with boundary. A simple choice would be given by two copies of $\R^d$ smoothly attached one to another through a compact corridor, tube or collar. As explained in \autoref{sec:intro}, in this work we focus on the scaling behavior of the constant in the Poincar\'{e} inequality for bilinear forms on dumbbell shaped subdomains of $\R^d$ satisfying \hyperlink{ConditionA}{Condition A}, see for example \hyperlink{figure1}{Figure 1}.

\begin{example}\label{ex:1}
Let us provide two domains satisfying \hyperlink{ConditionA}{Condition A}. The first one satisfies \eqref{eq:Condition_C}, the second one does not.  Set
 \begin{align*}
    &\Gamma = \{ x=(x_1,\tilde x) \in \R\times\R^{d-1} | \, |\tilde x| < 1\}\\ 
    & D=\{x\in\R^d | \, x_1 < -1 \} \cup \Gamma \cup \{x\in\R^d | \, x_1 > 1 \}.
 \end{align*}
Then $D$ satisfies \hyperlink{ConditionA}{Condition A} and \eqref{eq:Condition_C}. $D$ does not satisfy \eqref{eq:Condition_C} if we substitute $\Gamma$ by a non-convex uniform set without visibility through the corridor, e.g.~if
\begin{align*}
\Gamma=\{ x=(x_1,\tilde x) \in \R\times\R^{d-1} | \, |\tilde x-(2x_1^2-2,0,..,0)|< 1\}.
\end{align*}
In both cases, given $R>0$, we define $D_R$ by $D_R= D\cap B(0,R)= D^-_R \cup \Gamma \cup D^+_R$, where
\begin{align*}
\begin{split}
& D^-_R := \{x \in B(0,R) | x_1 < -1 \}, \quad  D^+_R = 
\{x \in B(0,R) | x_1 > 1 \}.
\end{split}
\end{align*}
\end{example}

\medskip

Our aim is to compare the scaling behavior of the Poincar\'{e} constant for local quadratic forms and nonlocal forms with visibility contraint. Let us first provide a Poincar\'{e} inequality in the local case. Note that we were not able to find a proof of this result although the result is stated at several places and seems to be well known. 

\begin{theorem}\label{thm:poincare-local}
Let $D$ be a domain satisfying \hyperlink{ConditionA}{Condition A}. Then there exists a constant $R_0>0$ such that for all $R\geq R_0$, $1\leq p\leq d$ and every $u\in W^{1,p}(D_R)$
\begin{align*}
\int_{D_R}|u(x)-u_{D_R}|^p \dx\lesssim 
\begin{cases}
R^d\int_{D_R}|\nabla u(x)|^p \dx,&d> p \\
R^p(\log R)^{p-1}\int_{D_R}|\nabla 
u(x)|^p \dx,&d=p. 
\end{cases}
\end{align*}
\end{theorem}

 \proof
Recall that $\Gamma^*=\Gamma\setminus (D^+\cup D^-)$. Let $a>0$, $L>0$ be such that there exist balls $B_a^+\subset D^+\cap\Gamma_L$, $B_a^-\subset D^-\cap\Gamma_L$ of radius $a$, where $\Gamma_L$ is a uniform subset of $D$ with diameter $L$ such that $\Gamma^* \cup B_a^+\cup B_a^-\subset \Gamma_L$. Then there is a collection of increasing sets $(C_{i})_{i\leq k}$ such that $C_0=B_a^+$, $C_k=D_R^+$ and
\begin{align*}
 |C_i|\asymp 2^{id}\text{ and }\diam(C_i)\asymp 2^i,\ i=0,1,...,k
\end{align*}
where $k\asymp \log R$. First, we compare average values of $u$ on sets $B_a^\pm$ and $\Gamma_L$, 
\begin{align}
  |u_{B_a^\pm}-u_{\Gamma^*}|^p &\lesssim |u_{B_a^\pm}-u_{\Gamma_L}|^p+|u_{{\Gamma_L}}-u_{\Gamma^*}|^p\nonumber \\
  &
\leq \frac{1}{|B_a^\pm|}\int_{B_a^\pm}|u(x)-u_{\Gamma_L}|^p 
\dx+\frac{1}{|{\Gamma^*}|}\int_{{\Gamma^*}}|u(x)-u_{\Gamma_L}|^p \dx\nonumber\\ 
  &\lesssim \int_{\Gamma_L}|u(x)-u_{\Gamma_L}|^p \dx\lesssim \int_{{\Gamma_L}}|\nabla u(z)|^p 
dz,\label{eq:means} 
\end{align}
where the last inequality follows by applying the classical Sobolev-Poincar\'{e} inequality on a uniform domain ${\Gamma_L}$, see for example \cite{Mar88}. Furthermore, 
\begin{align}
 \int_{D_R^\pm}|u(x)-u_{B_a^\pm}|^p \dx&\lesssim  
\int_{D_R^\pm}|u(x)-u_{D_R^\pm}|^p 
\dx+\int_{D_R^\pm}\left(\sum_{i=0}^k|u_{C_i}-u_{C_{i-1}}|\right)^p 
\dx\nonumber\\ 
 &\lesssim  R^p\int_{D_R^\pm}|\nabla u(x)|^p 
\dx+R^d\left(\sum_{i=0}^k|u_{C_i}-u_{C_{i-1}}|\right)^p.\label{eq:RBMproof1} 
\end{align}
By H\"{o}lder inequality and the classical Poincar\'{e} inequality, we have
\begin{align*}
 \sum_{i=0}^k|u_{C_i}-u_{C_{i-1}}|&=\sum_{i=0}^k\frac 
1{|C_{i-1}|}\int_{C_{i-1}}|u(x)-u_{C_{i}}|\dx\\
&\lesssim 
\sum_{i=0}^k2^{-id/p}\left(\int_{C_{i}}|u(x)-u_{C_{i}}|^p \dx\right)^{1/p}\\ 
 &\lesssim \sum_{i=0}^k2^{-id/p}2^{i}\left(\int_{C_{i}}|\nabla u(x)|^p 
\dx\right)^{1/p}\\
&\lesssim 
\left(\sum_{i=0}^k2^{-i(d-p)/p}\right)\left(\int_{D_{R}^\pm}|\nabla u(x)|^p 
\dx\right)^{1/p}\,.
 \end{align*}
 This implies for $p<d$
\begin{align*}
  \sum_{i=0}^k|u_{C_i}-u_{C_{i-1}}|\lesssim\left(\int_{D_{R}^+}|\nabla 
u(x)|^p \dx\right)^{1/p}             
\end{align*}
and for $d=p$
\begin{align*}
  \sum_{i=0}^k|u_{C_i}-u_{C_{i-1}}|\lesssim k^\frac 1 
q\left(\int_{D_{R}^+}|\nabla u(x)|^p \dx\right)^{1/p}\asymp \log(R)^\frac 
{p-1}p\left(\int_{D_{R}^+}|\nabla u(x)|^p \dx\right)^{1/p}. 
\end{align*}
Since
\begin{align}
 \int_{D_R}|u(x)-u_{D_R}|^p \dx & \leq 
\int_{D_R^+}|u(x)-u_{\Gamma^*}|^p \dx\nonumber\\
& \leq 
\int_{D_R^+}|u(x)-u_{B_a^+}|^p \dx+R^d|u_{B_a^+}-u_{{\Gamma^*}}|^p\nonumber\\
& \quad +\int_{{\Gamma^*}}
|u(x)-u_{{\Gamma^*}}|^p \dx+\int_{D_R^-}|u(x)-u_{B_a^-}|^p \dx\nonumber\\
& \quad +R^d|u_{B_a^-}-u_{{\Gamma^*}}|^p, \label{eq:RBMproof2}
 \end{align}
the proof of the theorem now follows from the calculation above together with 
\eqref{eq:means} and \eqref{eq:RBMproof1}. 
\qed  

\begin{remark}
With regard to sharpness of \autoref{thm:poincare-local} we present the following example. Assume that $D$ is the domain from \autoref{ex:1}. Given $R > 0$ sufficiently large, define a function $u:D_R \to \R$ by
\begin{align}\label{eq:sharpness}
u(x):=\begin{cases}
       -1,&x\in D_R^-\\
       x_1,& x\in \Gamma^*\\
       1,&x\in D_R^+.
       \end{cases}
 \end{align}
One easily checks $u_{D_R}=0$, $||u||^p_{L^p(D_R)}\asymp R^d$ 
and $\int_{D_R}|\nabla u(x)|^p\dx\asymp |\Gamma^*|$, which shows that the 
Poincar\'{e} constant is at least of order $R^d$.  For a more general domain, one can can construct an analogous example by taking a smooth function $u$ such that $u_{D_R}=0$, $u=\pm1$ on $D_R^\pm$ and $\nabla u$ is bounded on $\Gamma^*$.
\end{remark}
 
%The following lemma is needed below. Its proof is straightforward.  
%\begin{lemma}
%	Let $D$ be an open bounded set. Then the following Poincar\'{e} inequality holds true:
%	\begin{align*}
%	\int_{D}|u(x)-u_{D}|^p\dx&\lesssim c(d,s,p) \text{diam}(D)^{sp}\int_D\int_D\frac{|u(x)-u(y)|^p}{|x-y|^{d+ps}}\qquad (u\in L^p(D)).
%	\end{align*}
%\end{lemma}

Finally, we provide a proof of \autoref{thm:poincare-intro}. We formulate a more general result allowing for kernels satisfying conditions analogous to \eqref{eq:K2intro}, \eqref{eq:l1}, \eqref{eq:l2}, and \eqref{eq:ell1} in the $L^p$-setting. \autoref{thm:poincare-intro} then is just a corollary.

\begin{theorem}\label{thm:poincare-nonlocal}
Let $D$ be a domain satisfying \hyperlink{ConditionA}{Condition A}, where $\Gamma$ is a uniform domain. Let $p\geq 1$ and let the function $\ell:(0,\infty)\to(0,\infty)$ satisfy
\begin{align*}
&\int_0^{\infty}\left(r^{p-1}\wedge\frac 1 r\right)\ell(r)\d r<\infty,\\
&\lambda^{-\gamma}\lesssim \frac{\ell(\lambda r)}{\ell(r)}\lesssim \lambda^{-\delta}\qquad (\lambda\geq 1,\, r>0),
\end{align*}
for some constants $0 < \delta \leq \gamma<d$. Then there exists $R_0>0$ such that for every $R\geq R_0$ and $u\in  L^p(D_R)$, 
\begin{align*}
\int_{D_R}|u(x)-u_{D_R}|^p \dx\lesssim  
R^d\int_{D_R}\int_{D_{R,x}}|u(y)-u(x)|^p\frac{\ell(|x-y|)}{|x-y|^{d}}\dy \,
\dx. 
\end{align*}
If $D$ additionally satisfies condition \eqref{eq:Condition_C} and $\ell(R)\geq R^{-1}$ for $R\geq R_0$, then
\begin{align*}
\int_{D_R}|u(x)-u_{D_R}|^p \dx\lesssim R^{d-1}\ell(R)^{-1}\int_{D_R}\int_{D_{R,x}}|u(y)-u(x)|^p\frac{\ell(|x-y|)}{|x-y|^{d}}\dy \,
\dx. 
\end{align*}
\end{theorem}

\begin{remark}
There is an alternative way to formulate the result. If $\Gamma$ satisfies the stronger assumption \hyperlink{ConditionB}{Condition B}	 and $\ell:(0,\infty)\to(0,\infty)$ the weaker assumption 
\begin{align}\label{eq:p}
\lambda^{-\gamma}\lesssim \frac{\ell(\lambda r)}{\ell(r)}\lesssim \lambda^{d},\qquad (\lambda\geq 1,\,r>0)
\end{align}
for some $\gamma< d$ instead of \eqref{eq:l2}, then the assertion remains true. One would only need to work with a generalisation of \autoref{thm:comparability-weakly} in the $L^p$-setting, instead of a generalisation of \autoref{thm:comparable-intro}.
\end{remark}

\proof We apply the notation from the proof of \autoref{thm:poincare-local}. 
Let $C$ be a bounded set in $\R^d$ such that $|C|\gtrsim \diam(C)^d$. It is easy to see that the following Poincar\'{e} inequality holds,
\begin{align*}
 \int_{C}|u(x)-u_{C}|^p \dx\lesssim \ell(\diam(C))^{-1} 
\int_{C}\int_{C}|u(y)-u(x)|^p\frac{\ell(|x-y|)}{|x-y|^{d}}\dy \dx \,. 
\end{align*}
By applying the Poincar\'{e} inequality in the last line of \eqref{eq:means} we obtain 
 \begin{align}
  |u_{B_a^\pm}-u_{\Gamma^*}|^p &\lesssim 
\int_{\Gamma_L}\int_{\Gamma_L}|u(x)-u(y)|^p\frac{\ell(|x-y|)}{|x-y|^{d}}\dy \dx\nonumber\\
&\lesssim \int_{\Gamma_L}\int_{\Gamma_{L,x}} |u(x)-u(y)|^p\frac{\ell(|x-y|)}{|x-y|^{d}}\dy \dx.\label{eq:means4}
 \end{align}
where the second inequality follows from a straightforward generalisation of \autoref{thm:comparable-intro} in the $L^p$-setting. Similarly as in the proof of \autoref{thm:poincare-local} we obtain 
\begin{align}
 \int_{D_R^\pm}|u(x)-u_{B_a^\pm}|^p \dx&\lesssim  
\ell(R)^{-1}\int_{D_R^\pm}\int_{D_R^\pm}|u(x)-u(y)|^p\frac{\ell(|x-y|)}{
|x-y|^{d}}\dy \dx \nonumber\\
& \quad + R^d\left(\sum_{i=1}^k|u_{C_i}-u_{C_{i-1}}|\right)^p\nonumber 
\end{align}
and by \eqref{eq:p},
\begin{align*}
 \sum_{i=1}^k|u_{C_i}-u_{C_{i-1}}|&\lesssim 
\sum_{i=1}^k2^{-id/p}\ell(2^{i})^{-1/p}\left(\int_{C_{i}}\int_{C_{i}}
|u(x)-u(y)|^p\frac{\ell(|x-y|)}{|x-y|^{d}}\dy \dx\right)^{1/p}\\ 
 &\lesssim 
\left(\int_{D_{R}^+}\int_{D_{R}^+}|u(x)-u(y)|^p\frac{\ell(|x-y|)}{|x-y|^{
 d}}\dy \dx\right)^{1/p}. 
\end{align*}
These calculations together with \eqref{eq:means4} and \eqref{eq:RBMproof2} 
give the first inequality. 

Next, assume that $D$ additionally satisfies condition \eqref{eq:Condition_C}. Let $\Gamma_R^\pm:=D_R^\pm\cap 
\widetilde\Gamma$ and $K_R:=\Gamma_R^-\cup(\widetilde\Gamma\cap\Gamma^*)\cup\Gamma_R^+$ and note that these sets are convex. Similarly as above, it follows 
  \begin{align*}
\begin{split}
 \int_{D_R}|u(x)-u_{D_R}|^p\dx&\lesssim 
\int_{D_R^+}|u(x)-u_{\Gamma_R^+}|^p 
\dx+R^d|u_{\Gamma_R^+}-u_{K_R}|^p\\ 
 & \quad +\int_{\Gamma^*}|u(x)-u_{\Gamma^*}|^p \dx+|\Gamma^*|\cdot|u_{\Gamma^*}-u_{K_R}|^p \\
 & \quad +\int_{D_R^-}|u(x)-u_{\Gamma_R^-}|^p \dx+R^d|u_{\Gamma_R^-}-u_{K_R}|^p 
\end{split}
 \end{align*}
 and
 \begin{align*}
  |u_{\Gamma_R^\pm}-u_{K_R}|^p &\leq 
\frac{1}{|\Gamma_R^\pm|}\int_{\Gamma_R^\pm}|u(x)-u_{K_R}|^p \dx\nonumber\\
&\lesssim 
\frac{\ell(\diam(K_R))^{-1}}{ 
|\Gamma_R^\pm|}\int_{K_R}\int_{K_R}|u(x)-u(y)|^p\frac{\ell(|y-x|)}{|x-y|^{d}}\d 
y\, 
\dx\nonumber\\ 
  &\lesssim 
(\ell(R)R)^{-1}\int_{K_R}\int_{K_R}|u(x)-u(y)|^p\frac{\ell(|y-x|)}{|x-y|^{d}}\dy \, 
\dx.
 \end{align*}
As before we obtain 
\begin{align*}
 \int_{D_R^\pm}&|u(x)-u_{\Gamma_R^\pm}|^p \dx\lesssim 
\int_{D_R^\pm}|u(x)-u_{D_R^\pm}|^p\dx+R^d|u_{D_R^\pm}-u_{{
\Gamma_R^\pm}}|^p\nonumber\\ 
 &\lesssim  
\ell(R)^{-1}\int_{D_R^\pm}\int_{D_R^\pm}|u(x)-u(y)|^p\frac{\ell(|y-x|)}{
|x-y|^{d}}\dy \dx+\frac{R^d}{|\Gamma_R^\pm|}\int_{\Gamma_R^\pm}|u(x)-u_{
D_R^\pm}|^p \dx\nonumber\\ 
 &\lesssim  
(\ell(R)^{-1}+\ell(R)^{-1}R^{d-1})\int_{D_R^\pm}\int_{D_R^\pm}
|u(x)-u(y)|^p\frac{\ell(|y-x|)}{|x-y|^{d}}\dy \dx,\nonumber 
\end{align*}
and
\begin{align*}
\int_{\Gamma^*}|u(x)-u_{\Gamma^*}|^p \dx& \lesssim\int_{\Gamma^*}\int_{\Gamma^*_x}|u(x)-u(y)|^p\frac{\ell(|x-y|)}{|x-y|^{d}}\dy \dx \\
&\lesssim\int_{\Gamma^*}\int_{\Gamma^*_x}|u(x)-u(y)|^p\frac{\ell(|x-y|)}{|x-y|^{d}}\dy \dx,
\end{align*}
where the last line follows analogously as \eqref{eq:means4} by applying the generalisation of \autoref{thm:comparable-intro}. These inequalities together with 
 \begin{align*}
|u_{\Gamma^*}-u_{K_R}|^p  &\lesssim |u_{\Gamma^*}-u_{\widetilde\Gamma}|^p+|u_{\widetilde\Gamma}-u_{K_R}|^p \\
&\leq \frac{1}{|\widetilde\Gamma|}\int_{\widetilde\Gamma}|u(x)-u_{\Gamma^*}|^p \dx+
\frac{1}{|\widetilde\Gamma|}\int_{\widetilde\Gamma}|u(x)-u_{K_R}|^p \dx\nonumber\\
&\lesssim \int_{\Gamma^*}|u(x)-u_{\Gamma^*}|^p \dx+
\int_{K_R}|u(x)-u_{K_R}|^p \dx\nonumber\\
&\lesssim \int_{\Gamma^*}\int_{\Gamma^*}|u(x)-u(y)|^p\frac{\ell(|y-x|)}{|x-y|^{d}}\d 
y\, 
\dx\nonumber\\
&\quad + \ell(\diam(K_R))^{-1}\int_{K_R}\int_{K_R}|u(x)-u(y)|^p\frac{\ell(|y-x|)}{|x-y|^{d}}\d 
y\, 
\dx\nonumber\\
&\lesssim \ell(R)^{-1}\int_{D_R}\int_{D_{R,x}}|u(x)-u(y)|^p\frac{\ell(|y-x|)}{|x-y|^{d}}\d 
y\, 
\dx\nonumber
 \end{align*}
finish the proof.
\qed  

\proofof \autoref{thm:poincare-intro}: The theorem follows from \autoref{thm:poincare-nonlocal} by its application to the function $\ell(r)=r^{-sp}$ for $s\in(0,1)$. 
\qed

\begin{remark}
 In order to show that these inequalities are sharp for $s\in(0,1)$, we consider $D$ as the domain from \autoref{ex:1}, which does not satisfy condition \eqref{eq:Condition_C}, and $u$ as in \eqref{eq:sharpness}. Let $x_0 \in \Gamma^*$ and $R \geq 1$. Then 
\begin{align}
\int\limits_{D_R}\int\limits_{D_{R,x}}\frac{|u(x)-u(y)|^p}{|x-y|^{d+sp}}\dy \dx &= \int\limits_{\Gamma^*}\int\limits_{\Gamma^*}\frac{|x_1-y_1|^p}{|x-y|^{d+sp}}\dy \dx 
+ 2\int\limits_{\Gamma^* }\int\limits_{D_{R,x}^\pm}\frac{|x_1\mp1|^p}{|x-y|^{d+sp}}\dy \dx\nonumber\\  
&\lesssim\int_{\Gamma^*}\int_{\R^d}\frac{(1\wedge 
|z|^p)}{|z|^{d+sp}}dz\dx\lesssim \frac{|\Gamma^*|}{s(1-s)}.\label{eq:sharp1}
\end{align}
Next assume $D$ is the domain from \autoref{ex:1} satisfying the condition \eqref{eq:Condition_C}. Then
\begin{align*}
&\int\limits_{D_R}\int\limits_{D_{R,x}}\frac{|u(x)-u(y)|^p}{|x-y|^{d+sp}}\dy \dx \\
&\asymp \int\limits_{\Gamma^*}\int\limits_{\Gamma^*}\frac{|x_1-y_1|^p}{|x-y|^{d+sp}}\dy \dx + \int\limits_{\Gamma^* }\int\limits_{D_{R,x}^\pm}\frac{|x_1\mp 1|^p}{|x-y|^{d+sp}}\dy \dx + \int_{D_{R}^-}\int_{D_{R,x}^+}\frac{1}{|x-y|^{d+sp}}\dy \dx \\
&=:I_1+I_2+I_3. 
\end{align*}
Analogously as in \eqref{eq:sharp1},
\begin{align*}
I_1+I_2\lesssim\int_{\Gamma^*}\int_{\R^d}\frac{(1\wedge 
|z|^p)}{|z|^{d+sp}}dz\dx\leq \frac{|\Gamma^*|}{s(1-s)}.
\end{align*}
For the remainder term, define
\begin{align*}
 K^-:=\bigcup\limits_{y\in D_R^+ }D_{R,y}^+\cap 
D_R^-=\{x=(x_1,\tilde x)\in D_R^-|\, -R<x_1<-1,\,|\tilde x|\leq c_1x_1\},
\end{align*}
for some constant $c_1>0$ which is independent of $R$. Note that for $x\in K^-$ 
\begin{align*}
D_{R,x}^+\subset V_{\alpha(x_1)}(x, x_1, 2R)
\end{align*}
where $\tan\frac{\alpha(x_1)}2\asymp \frac 1{x_1}$. Therefore,
\begin{align*}
I_3&\asymp \int_{K^-}\int_{D_{R,x}^+}\frac{1}{|x-y|^{d+sp}}\dy \dx\lesssim 
 \int_{K^-}\int_{V_{\alpha(x_1)}(x,x_1,2R)}\frac{1}{|x-y|^{d+sp}}\dy \dx\\ 
&\lesssim \int_{K^-}\frac{1}{x_1^{d-1}}\int_{x_1}^{2R}\frac{1}{r^{1+sp}}\d r\d 
x 
\lesssim \int_1^{R}\int_{|\tilde x|<c_2x_1}\frac{1}{x_1^{d+sp-1}}\d\tilde 
x\dx_1\\
&\lesssim 
\begin{cases}
\frac{ R^{1-sp}}{1-sp},&s\in(0,1/p)\\
%\log(R),&s= 1/p \\ 
\frac{ 1}{1-sp},&s\in(1/p,1)
\end{cases}
\end{align*}
Since $u_{D_R}=0$, $||u||_{L^p(D_R)}^p\asymp R^d$ and
\begin{align*}
\int_{D_R}\int_{D_{R,x}}\frac{|u(y)-u(x)|^2}{|x-y|^{d+2s}}\dy \,
\dx\lesssim 
\begin{cases}
R^{1-sp},&s\in(0,1/p)\\
%\log(R),&s= 1/p \\
1,&s\in(1/p,1)
\end{cases},
\end{align*}
we obtain that the dependence on $R$ of the constant in \autoref{thm:poincare-intro} is sharp for $s\neq \frac{1}{p}$. 
\end{remark}

%\bibliographystyle{plain}
%\bibliography{geometry_constrained}

\def\cprime{$'$}

% \listoftodos

\end{document}